\documentclass[11pt]{article}
\usepackage{amsmath,amssymb}

\usepackage{fullpage}
\usepackage{color}
\usepackage{amsmath}

\newcommand{\qed}{ $\blacksquare$}

           \newcommand\om{\omega}
\newcommand{\var}{\varepsilon}

\newcommand{\pa}{\partial}

\def\eps{\varepsilon}

\Large

\newtheorem{theo}{Theorem}

\newtheorem{lem}[theo]{Lemma}
\newtheorem{rem}[theo]{Remark}
\newtheorem{defi}[theo]{Definition}
\newtheorem{cor}[theo]{Corollary}
\numberwithin{equation}{section}

\begin{document}

\title{On a constrained 2-D Navier-Stokes Equation}

\author{Caglioti, Pulvirenti \footnote{Dipartimento di
Matematica, Universit\`a di Roma `La Sapienza', P.le Aldo Moro 2,
\hfill\break 00185 Roma, Italy. E-mail: {\tt caglioti@mat.uniroma1.it, pulvirenti@mat.uniroma1.it} },
Rousset\footnote{CNRS, Laboratoire Deudonn\'e,   Universit\'e de Nice, Parc Valrose, 06108 Nice
cedex 2, France.  E-mail: {\tt frousset@unice.fr}.}
}

\date{}

\maketitle

\begin{abstract}
The planar Navier-Stokes equation exhibits, in absence of external
forces, a trivial  asymptotics in time. Nevertheless the appearence of
coherent structures suggests non-trivial intermediate asymptotics
which should be explained in terms of the equation itself.
Motivated  by  the separation of the 
different time scales  observed in  the  dynamics of the Navier-Stokes
equation, we study  the well-posedness and asymptotic behaviour of a constrained equation which  neglects the
 variation of the  energy and  moment of inertia.

\end{abstract}

\section{Introduction}
Consider the  two-dimensional Euler equation in vorticity form 
\begin{equation}
\label{euler}
(\pa_t+u\cdot \nabla)\om(x,t)=0, \quad x \in \mathbb{R}^2
\end{equation}
where 
 the divergence free velocity field $u$ is given by 
  $u= \nabla^\perp \psi$, $\psi = - \Delta^{-1} \omega$.
  Explicitely, we can write: 
\begin{equation}
\label{K}
u=K*\om, \qquad K(x)= \nabla^{\perp}g =-\frac 1{2\pi} \frac {x^{\perp}}{|x|^2}, 
\quad g(x)= - { 1 \over 2\pi} \log |x|.
\end{equation}
The rigorous justification of the formation of coherent structures
 in  two-dimensional fluid-dynamics, which is 
 observed  in real and numerical experiments (see  e.g. \cite{MJnum}) remains
 a widely open problem.
 An  attempt to justify the appearance of these coherent structure
    is due to Onsager \cite{O}, see also \cite{LP}, \cite{MJ}, and  \cite{ES} for a recent review. The main idea
     is to  replace the incompressible  Euler equation  by the system of
     $N$  point vortices  and to study  the Statistical Mechanics
      of these point vortices.  In the mean field  limit $N\rightarrow +\infty$, 
      the Gibbs measure associated to the point vortices concentrates to
         some special stationary solutions of the Euler equation
         (called mean field solutions), we refer to   \cite{CLMP1},
         \cite{CLMP2},\cite{Kie},\cite{KieL} for the rigorous justification.
These states are under the form: 
\begin{equation}
\label{MF}
\om=\frac {e^{b \psi + a \frac {|x|^2}2}} {Z}, \quad Z=\int_{\mathbb{R}^2} e^{b  \psi+ a  \frac {|x|^2}2}\, dx.
\end{equation}
In this  last expression, $Z$ 
is a normalization factor  to have $\int \om = 1$   and
 $b$  real and $a <0$ are  parameters.
 From a mathematical point of view  this equation enters in 
 the framework  of a  general class of   nonlinear elliptic equations  given by
\begin{equation}
\label{MFeq}
-\Delta \psi =\frac {e^{ b  \psi + a |x|^2}}{Z}
\end{equation}
which has been studied in \cite{CLMP1}, \cite{CLMP2}.
Nevertheless, there is no justification of the fact that,   among the infinite number of stationary stable
solutions of the Euler equation,   the mean field solutions
play  indeed  a special role in the dynamics.

Another justification of \eqref{MF} could come from the intermediate asymptotic
 behaviour of the two-dimensional  Navier-Stokes equation which,  in vorticity form, reads
 
 \begin{equation}
\label{nstokes}
(\pa_t+u\cdot \nabla)\om(x,t)=\nu \Delta \om (x,t), \quad x \in \mathbb{R}^2, 
\end{equation}
where  $\nu>0$ is the viscosity coefficient. Indeed, 
 due to the dissipation term in the right hand side of eq.n \eqref{nstokes}, the
asymptotic behaviour of the solutions is trivial, namely, when 
 $t \rightarrow + \infty$,  $\om (x,t)\to 0$
pointwise and in the $L^p$ sense,  for $p>1$.   Consequently,
  the states \eqref{MF} could play a part only in the intermediate behaviour
 of the equation before the dissipation scale. To give a more quantitative
  description of this idea, it is usefull to recall that
   the  solutions of \eqref{MF}  can be  studied through a variational
   principle: the radial solutions  are obtained
   as minimizers of
    the Boltzmann entropy (this is proven in the Appendix)
     $$S(\om)= \int_{\mathbb{R}^2} \om \log \om \, dx$$
      under the constraints
       $$ E(\om)= E, \quad M(\om) =  0, \quad  I(\om)= I, \quad \om \geq 0,  \quad \int_{\mathbb{R}^2} \om\, dx =1,
       $$
       for some fixed $E$ and $I$,
       where  the energy $E(\om)$, the center of  vorticity  $M(\om)$ and the moment
        of inertia $I(\om)$ are  respectively  given by
      $$ E(\om)= { 1 \over 2} \int_{\mathbb{R}^2} \psi\, \om \, dx, \quad
      M(\om) = \int_{\mathbb{R}^2} x \om \, dx, \quad
       I(\om)={{ 1 \over 2}} \int_{\mathbb{R}^2} |x-M|^2 \om.$$
       Note that it is always possible to choose the coordinates so that
        $M(\om)=0$.
       
       It is thus interesting to study how these quantities, which are
        conserved by the Euler equation \eqref{euler}, evolve   under   the Navier-Stokes 
         flow. At first, it is well known
         that the Navier-Stokes equation preserves the nonnegativity
       and that $\int \omega$ and $M$ are conserved.
        Consequently, throughout this paper, we will focus
         on non-negative solutions which are normalized  such that $\int \om =1$ and            $M(\om)=0$.
        Next, we also observe that 
   $$ \dot I(\om)=2 \nu, \quad   \dot E=- \nu \int_{\mathbb{R}^2} \om^2, \quad    \dot S(\om)  =    - \nu \int_{\mathbb{R}^2} {|\nabla \omega |^2 \over \omega}.$$
   It is easy to see that $I$  can be considered as constant for times
    $t << \nu^{-1}$.  Moreover,  it is likely that in certain cases  
    (see for instance \cite{MJnum})
    the energy dissipation rate is much smaller than the entropy dissipation rate.
    %also seems that the energy dissipation rate 
    % very often  is much smaller than the entropy dissipation rate.        
 Coming back to an attempt to justify eq. \eqref{MFeq}
in terms of the Navier-Stokes evolution, the first
naive remark is that, if the energy and  the moment of
inertia are assumed to  vary  on a long time scale,  
they can be considered as constant in a first approximation. On such
a time scale,  the motion should be governed by a master equation,
which modifies the Navier-Stokes equation leaving constant both energy
and moment  of inertia, but retaining all the other features of the
Navier-Stokes dynamics. Therefore such a master equation,
dissipating the entropy at constant energy and moment of inertia,
would lead to the solution to eq.n \eqref{MFeq} as $t \to \infty$.
%Such a prescription is, however, ambiguous, as we shall see later on by
%discussing a simple example.
 By using  a recent geometric  gradient flow 
characterization of the Navier-Stokes equation (see \cite {V} for example) connected with
  the mass
transport problem and the associated differential calculus introduced
 in \cite{Otto} (see also \cite {AGS}) we have derived such an equation
  in \cite{CPR}. 
 In this framework
 the Navier-Stokes equation can be written as a differential equation
  for a  vector field which is 
  the sum of a dissipative part,  which is the    gradient flow 
  of the entropy,  and a conservative part,  
corresponding to the Euler equation,
which is the  orthogonal gradient of the energy.
  The equation we  were looking for was  obtained  by keeping only 
  the orthogonal projection of  the vector field in the tangent space of
  the  manifold $I=$const and $E=$const.  Due to the  nature of this decomposition, this procedure
  modifies the dissipative part while leaving  invariant the conservative part.
  Thus we have  found the equation:
\begin{eqnarray}
\label{nscons}
\pa_t \om +u\cdot \nabla \om & = &\nu \,  \hbox {div}(\nabla \om -
b\,  \om \nabla\psi-a\,  \om x)   \\
\nonumber & =&  \nu\,  \hbox {div}\bigl[\om \nabla (\log \om -
b\psi-a \frac {|x|^2}2)\bigr],
\end{eqnarray}
where the Lagrange multipliers $b$ and $a$ are given by 
\begin{equation}
\label{ba}
b = b(\om) =\frac {2I \int \om^2 +2V}{2I \int \om |\nabla\psi|^2 -V^2}; \quad
a = a(\om) =-\frac {2 \int \om |\nabla\psi|^2+V \int \om^2}{2I \int \om
|\nabla\psi|^2 -V^2},
\end{equation}
and
\begin{equation}
\label{V}
V=\int \om x \cdot \nabla \psi=\int_{\mathbb{R}^2} \! \int_{\mathbb{R}^2}  \om (x, t) \om (y, t )\, 
x \cdot \nabla g(x-y)\, dx dy =-\frac 1 {4\pi}.
\end{equation}

A way to validate this approach, is to  test it in the following simpler
 and well understood case.
 
A special self-similar  solution to eq.n \eqref{nstokes}  is 
 the so called Oseen vortex:
$$
 \omega (x,t)=\frac 1 {4 \pi  \nu (t+1)} e^{-\frac {|x|^2}{4\nu (t+1)}}.
$$
 Note that this is also a solution to the heat equation. 
 It was shown  by Gallay and Wayne \cite{GW} that this solution describes the long time asymptotic of the Navier-Stokes  equation in $L^1$.
 Indeed, with the  change of variables
$$
 \xi =  \frac {x}{\sqrt {1+t}}; \qquad \tau =   \log (1+t),  \qquad 
 \om (x,t)=(1+t)^{-1} \, w (\xi, \tau),
$$
the Navier-Stokes equation in the new variables reads :
\begin{equation}
\label{nsresc}
\partial_{\tau} w  + v \cdot \nabla_{\xi} w = 
 \nu  \Delta_{\xi} w +\nu\,\nabla_{\xi} \cdot \Bigl(\frac 12 \xi w \Bigr). 
 \end{equation}
 It is possible to show that $w \to W$ 
 in $L^1$ as $\tau \to \infty$,
where $W(\xi)$ 
is the rescaled Oseen vortex. As a consequence the Oseen vortex can be
thought as characterizing an intermediate asymptotics  for times  $\nu t <<1$.

This analysis enters perfectly in the context of
the projected gradient flows. Indeed imposing the
 constance of $I$  in the Navier-Stokes equation we find
\begin{equation}
\partial_{t} \omega  + u \cdot \nabla \omega =  \nu 
 \Delta \omega +\nu  \frac 1I \nabla \cdot\Bigl(\omega x\Bigr),  
\end{equation}
that is eq.n \eqref{nsresc} for $I=2$, as well as the rescaled Oseen vortex
$W(\xi)$   is  a Mean Field
solution  for $\beta=0$. 
This suggest to impose also the constance of $E$ in the attempt of outlining what happens before the occurrence of the Oseen vortex.  
 Indeed  one could argue that $I$ is more robust than $E$ in
many interesting physical situations. If so eq.n \eqref{nscons} should be more
appropriate on the time scale when $E$ is practically constant, while eq.n
\eqref{nsresc} should describe the fluid when $E$ start to be dissipated at
constant $I$. 

The aim of this paper is the mathematical study of  equations
 \eqref{nscons}, \eqref{ba}. For more details  on the derivation of 
   these equations   and of the physical motivations, we refer to \cite{CPR}.
    As explained in \cite{CPR}, the procedure of  constraining a diffusion
     equation is highly non unique. Nevertheless, an interesting
      feature   of  eq.n \eqref{nscons} is that it can be obtained by many
       different methods. As already explained, it appears naturally
        by using the  geometric structure of the Navier-Stokes equation.
        Moreover, 
   it was  noticed in \cite{CPR} that  eq.n \eqref{nscons} can also be obtained 
    by constraining the stochastic vortex dynamics
which is a finite dimensional approximation to the Navier-Stokes equation (see \cite {MP3},
 \cite {MPbook2} \cite {Osada}...).
Indeed (at a formal level), it is shown  that the stochastic process for a system of $N$
 stochastic  vortices, once constrained on the $E=$const manifold, produces, in the mean-field limit, exactly eq.n \eqref{nscons} which turns out to be compatible with this particle approximation.
 
   Equation \eqref{nscons} however is not new. It was previously derived by Chavanis in \cite{Ch1} and \cite {Ch2} following a  completely different approach based on the kinetic theory of  (deterministic) point vortices.

   We finally note that the projection on the manifold $I=$ const according to the gradient notion used here, has been considered by Carlen and Gangbo in  \cite{CG} for a different class of equations.

 A different approach to the one of Onsager \cite{O} in order  to understand the coherent structures arising in $2D$ flows
 was proposed by Robert and Sommeria \cite{RS} and Miller \cite{Miller}.
  The equilibrium solutions  are obtained 
  by using the maximum entropy principle  over a state space formed
   by a selected family of possible values of the vorticity. Note that this procedure preserves
    all the Euler invariants so that, as far as the equilibrium  is concerned, 
     the Robert-Sommeria-Miller  theory  is quite different 
             from  the  approaches,   as the present one,   based on the mean field equation.    As regards the dynamics, a class of master equations leading to such 
        equilibria
    has been introduced in  \cite{RS2}.  One of them,   which has some  formal
     similarities with our model,  has been
    systematically investigated from a  mathematical point of view in \cite{MR}.
     We remark that such an equation   exhibits  a maximum principle  leading to
      useful a priori estimates, for instance
      the $L^\infty$ norm of the vorticity is uniformly bounded  while,  in our case, we do not have 
      such a priori control.
      
  Our aim is to establish global existence results for \eqref{nscons} and
   to study the asymptotic behaviour of  global solutions. There are two main difficulties.
    The first one is that  
 eq.n \eqref{nscons} makes no sense whenever
the denominator in the definition of $a$ and $b$ (see eq.n \eqref{ba}) vanishes.
Note that by the Cauchy-Schwarz inequality we have 
$$
V^2=\Big(\int \om x \cdot \nabla \psi\Big)^2 \leq 2I  \int \om | \nabla \psi|^2
$$
and thus this denominator is always non-negative.
 Nevertheless, 
 when $\om \nabla
\psi$ and  $\om x$  are collinear, it vanishes. 
This happens for the one-dimensional family
of circular vortex patches:
$$
\omega =\frac 1 {\pi R^2} \chi_{B(0,R)}  
$$
where $\chi_{B(0,R)}  $ is the characteristic function of $B(0,R)$, the
disk of center $0$ and radius $R$. Indeed, we have:
$$
\omega \nabla \psi=-\frac {\om x }{2\pi R^2} \chi_{B(0,R)}.
$$
The other difficulty is that   $b$ is well-defined  if  $\om \in L^2$
 but  there is no a priori estimates available for the
 $L^2$ norm of the vorticity. The only a priori  information we have at our disposal are
  that $E$ and $I$ are  conserved (note that in this setting $E$ is not very usefull
   since the energy has no sign) and that the entropy $S$  decays. Indeed, we
    formally have 
 \begin{eqnarray*}
 \nonumber\frac {dS(\om)}{dt}  & = & - \nu \int \om \nabla \log \om  \cdot
 \nabla \big( \log \om - b \psi - a { |x |^2 \over 2 }  \big) \\
 &=& -\nu\, \int \om \big|\nabla\big( \log \om - a \frac {|x|^2}2 -b \psi \big)\big|^2.
\end{eqnarray*}
This identity can be checked by direct computation and it is obvious by using
 the geometric interpretation of the equation (see \cite{CPR}).
  Note that this identity is also usefull to guess that  asymptotic states
   should be given by the mean field solutions \eqref{MF} since the entropy
    dissipation  vanishes precisely on these states.
  
  Because of these difficulties, we  will be  able to get global
  existence results only for data sufficiently close to a mean field solution.
  Nevertheless,  we point out that  our  smallness constraint is independent
   of the viscosity  parameter $\nu$. It remains an open problem
     to establish if eq.n \eqref{nscons} can produce a singularity in a finite time
    and in particular if the  $L^2$ norm of $\om$  can blow up. Note that if we 
      consider \eqref{nscons} without \eqref{ba} i.e we consider the equation
       with some given parameters $a<0$ and $b$ fixed,  it is easy to 
        establish the existence of solution which blow-up.  Indeed,  since
         the inertial term does not play any part in the estimates, 
          all the result established in \cite{BDP} for the Keller-Segel equation
           remains true for this equation. In particular we have that 
            the evolution of $I$ which is nonnegative is given by
         $$ \dot I= a I + ( 2 - { b \over 4 \pi })$$
         and hence, if $b$ is larger than $8 \pi$, the solution must blow up in finite time.
       It would be very interesting to know if there is a nonlinear stabilization
        for \eqref{nscons}, \eqref{ba}.
  
 The paper is organized as follows:   
the global  existence proof is presented in Sect. 3 after that some preliminary steps are discussed in Sect. 2. Sect. 4 is devoted to the proof of 
 $L^p$ estimates necessary both for the existence part and the asymptotic behavior discussed in Sect. 5. Finally, the Appendix is devoted to the study of equation
  \eqref{MF} and   the connected variational principles. The main ideas
   are in \cite{CLMP1}, \cite{CLMP2}, but  the adaptation of these results
    to the  $\mathbb{R}^2$ case requires some care.

 \section{Preliminaries}

Let us introduce the  submanifold of probability densities
\begin{equation}
{\cal M}(E,I) =\Big\{\om, \quad \om \geq 0, \, \int \om =1, \, E(\om)=E,\,  I(\om)=I\Big\}
\end{equation}
for some fixed $E$ and $I>0$.
 Next,  by using Theorems \ref{theoMFE}, \ref{MFEmicro} in the Appendix, 
  we  denote the unique minimizer of the entropy functional $S(\om)$ on 
  $\mathcal{M}(E,I)$ by $\om_{MF}$. Note that $\om_{MF}$
   is  a  (radial) solution to eq.n  \eqref{MF} with parameters
    $b$ and  $a$ that we denote by  $b_{MF}$
  and   $a_{MF}$ respectively. 
\bigskip
\subsection{A stability property}
We  first  establish a crucial  result  asserting
the continuity of the $L^1$ norm with respect to variation of the
entropy, in a neighbourhood  of $\om_{MF}$, in the manifold ${\cal M}(E,I)$.

\begin{theo}
 \label{theoL1cont}
  For any $\var >0$ there exists $\delta >0$ such that, for all $\om  \in {\cal M}(E,I)$ for which
\begin{equation}
S(\om)-S(\om_{MF})\leq \delta,
\end{equation}
 then 
\begin{equation}
\|\om-\om_{MF} \|_{L^1}\leq \var.
\end{equation}
\end{theo}

%An interesting  by-product  of this result is the following: 
 %it was proven  in \cite{CLMP1} that $\om_{MF}$ is a (Lyapounov) stable
 % solution of the Euler equation \eqref{euler}. The argument uses  the fact that 
 %   it is a radial  stationary solution and thus that the result of \cite{MP2} which
  %    uses all the conserved quantities of the Euler equation applies.
 % Theorem \ref{theoL1cont} provides an alternative proof which uses only
  % the conservation of  $S$, $E$  and $I$ which are the natural quantities involved
   %in the problem.
      
  We remark that Thm 1 provides a proof of the Lyapounov stability of $\om_{MF}$  with respect to the Euler flow (by virtue of the time  invariance of $ S(\om )$ ). The stability of  $\om_{MF}$ was also proved in 
\cite{CLMP1}  by using the arguments of  \cite{MP2} in which
 all the conserved quantitities of the Euler equation are used.

{\bf Proof of Theorem \ref{theoL1cont}.} Assume, by contradiction, the existence of a sequence
$\om_n \in {\cal M}(E,I)$ and $\var >0$ such that
\begin{equation}
\label{hyp}
\lim_n S(\om_n)=S(\om_{MF}) \quad  \hbox {and} \quad
\|\om_n-\om_{MF} \|_{L^1}\geq \var.
\end{equation}

Thanks to the entropy bound, we can find a probability distribution $\om \in L^1$
such that, up to the extraction of a subsequence, 
$\lim_n \om_n=\om$ in the sense of the weak convergence of measures. Next, 
we also have that 
$$
\lim_n E(\om_n)=E(\om)=E
$$
(see the proof of  \eqref{En} in the Appendix ) and that 
$$
\lim_n I(\om_n)=I\geq I(\om), \quad \lim_n S(\om_n) \geq S(\om)
$$
by convexity.
Let now $a_{MF},\, b_{MF}$ be the multipliers  for which $\om_{MF}$ solves eq.n
(1.16) with those  values of parameters, and
 $ F_{(a_{MF}, b_{MF})}$ the free nergy functional (see \eqref{freeab} for the definition). 
  We  get, since  $a_{MF}<0,$ that 
\begin{eqnarray*}
F_{(a_{MF}, b_{MF})}(\om)  & \leq &  \lim_n \Big( S(\om_n)-b_{MF} E(\om_n)-a_{MF} I(\om_n)\Big) \\
&  = &( S(\om_{MF})-b_{MF} E(\om_{MF})-a_{MF} I(\om_{MF})) \\
& = & F_{(a_{MF}, b_{MF})}(\om_{MF}).
\end{eqnarray*}
Since $\om_{MF}$  is the unique minimizer of $F_{(b_{MF},a_{MF})}$
 (see Theorem \ref{theoMFE}),
  it follows that $\om=\om_{MF}$. As a consequence, we also get that 
\begin{equation}
\label{limI}
\lim_n I(\om_n)=I(\om) = I(\om_{MF}).
\end{equation}
Finally, we consider the relative entropy 
\begin{eqnarray*}
S(\om_n|\om_{MF}) & = & \int \om_n \log (\frac {\om_n }{\om_{MF}}) \\
& = & 
S(\om_n)-S(\om_{MF})+b_{MF} \int \Big(\om_{MF}-\om_n\Big)  \psi_{MF}+
 a_{MF}\Big( I(\om_{MF}) - I(\om) \Big).
\end{eqnarray*}
Now, we observe that   $S(\om_{n})-S(\om)$ goes to zero thanks to  \eqref{hyp}
 and that
$$ b_{MF} \int \Big(\om_{MF}-\om_n\Big)  \psi_{MF}+
 a_{MF}\Big( I(\om_{MF}) - I(\om) \Big)$$
  also goes to zero by weak convergence and \eqref{limI}. Consequently,
   the relative entropy  $S(\om_n|\om_{MF})$ goes to zero.
Thus we conclude by using the Csiszar-Kullback inequality
$$
\|\om_n-\om_{MF} \|_{L^1}^2 \leq  2  S(\om_n|\om_{MF})
\to 0
$$
which yields the desired contradiction.   \qed

\bigskip

\subsection{Properties of the coefficients $a(\om), b(\om). $}
Our next step is the study of the properties of $a$ and $b$.
We shall use the notation
\begin{equation}
\label{babis}
b(\om)=\frac {2I \int \om^2 +2V}{D(\om)},  \quad
a(\om) =-\frac {2 \int \om |\nabla\psi|^2+V \int \om^2}{D(\om) },
\end{equation}
where
\begin{equation}
D (\om)= 2I \int \om |\nabla\psi|^2 -V^2, \quad  V=-\frac 1{4\pi}.
\end{equation}
As we have seen in the introduction, one of the main difficulty is
 that the denominator $D(\om)$ may vanish for a vortex patch.

 Before stating the result,
we shall recall a usefull set of inequalities which will be used throughout the paper 
(see e.g.  \cite{Stein}...).
   \begin{lem}[Useful inequalities in $\Bbb R^2$] \hspace{1cm}
   \newline 
   \vspace{-0.5cm}
   \begin{itemize}
   \item[i)]{\bf Sobolev-Gagliardo-Niremberg : }
   The following inequalities hold for some $C>0$:
   \begin{eqnarray}
   \label{sob1}
 & &   ||\omega ||_{L^2} \leq C ||	\nabla \omega ||_{L^1} \\
    \label{sob2}
 & &   ||\omega||_{L^2}^2 \leq C \, ||\omega||_{L^1}\, ||\nabla \omega ||_{L^2}.
 \end{eqnarray}
 \item[ii)]{\bf Biot et Savart law : } 
  Let $u = K* \omega$ with $K$ defined by \eqref{K}, then we have 
  \begin{eqnarray}
& &  \label{vomega}
 1 < p < 2, \ 2<q< + \infty, \  \frac{1}{q} = \frac{1}{p} - \frac{1}{2}, \
 \vert \vert u \vert \vert _{L^q} \leq C \vert  \vert \omega   \vert \vert_{L^p} ,   \\
 & & \label{vinfty}
 1 \leq p <2, \ 2< q\leq \infty , \ \frac{1}{2} = \frac{\alpha}{p} + 
\frac{1- \alpha}{q} , \
 \vert \vert u \vert \vert_{L^\infty} \leq C \vert \vert \omega \vert 
\vert_{L^p}^\alpha 
  \, \vert \vert \omega \vert \vert_{L^q}^{1- \alpha},\\
  & & \label{nablav} 1<p < +\infty, \ \vert \vert \nabla  u \vert
  \vert_{L^p} \leq C \vert  \vert \omega \vert \vert_{L^p}.
  \end{eqnarray}
 \item[iii)]{\bf Interpolation in $L^p$ spaces:}
 \begin{equation}
 \label{pq}
1\leq  p<r<q\leq +\infty ,
 \quad ||\om||_{L^r} \leq 2 ||\om ||_{L^p}^\alpha\, || \om ||_{L^q}^{ 1- \alpha}, 
\quad \alpha= {p \over r} { 1- { r \over q } \over 1 -{ p \over q } }
 \end{equation}
 \end{itemize}
 \end{lem}

 We shall prove the following result: 
\begin{theo}
 \label{theoL1cont of a,b}
 Suppose $\om \in L^p \cap {\cal M}(E,I) $ for some $p\geq 2$. 
 \begin{enumerate}
\item We have:  
\begin{equation}
\label{bMF}
b(\om_{MF}) = b_{MF}, \quad a(\om_{MF}) = a_{MF}
\end{equation}
 and $D(\om_{MF})>0$.
 \item $D(\om)$, $a(\om)$ and $b(\om)$ are continuous on $L^2$
 \item Assume that $p>2$, then if 
  $S(\om) - S(\om_{MF})$ is sufficiently  small, we have
\begin{equation}
\label{DomMF}
|D(\om) - D(\om_{MF})|+ |a(\om)-a(\om_{MF}) |+|b(\om)-b(\om_{MF}) | \leq C \| \om-\om_{MF} \|_{L^1}^{\alpha},
\end{equation}
 for some  $\alpha <1$. Here $C$ depends on $E, I$ and $\|\om\|_{L^p} $.
 \end{enumerate}

\end{theo}

{\bf Proof of Theorem \ref{theoL1cont of a,b}.}

We first prove $1.$
Since  $\om_{MF}$ is  a solution of the mean field equation \eqref{MF}, 
 we have
 $$  \nabla \om_{MF} = b_{MF} \om_{MF} \nabla \psi_{MF} + a_{MF} \om_{MF} x.$$
 Taking  the scalar product of this equation
  by $\nabla \psi_{MF}$ and $x$   we find
$$ b_{MF} \int \om_{MF} |\nabla \psi_{MF}|^2 + a_{MF}V = \int \om_{MF}^2, \quad  b_{MF}V + 2a_{MF} I = -2.$$ 
The resolution of this two by two linear system precisely gives that
$$
 b(\om_{MF}) = b_{MF}, \quad a(\om_{MF}) = a_{MF}.
$$
Consequently,   since the numerator in the definition of $b(\om_{MF})$ and
 $a(\om_{MF})$ is finite, we  find that 
 \begin{equation}
 \label{posD}
 D(\om_{MF}) >0.
 \end{equation}
Next, we shall estimate  the  differences
\begin{eqnarray}
& & \label{dif1}
\int \om^2 -\int \om_{MF}^2 , \\  
& &  \label{dif2}
\int \omega |\nabla
\psi|^2-\int \omega_{MF} |\nabla \psi_{MF}|^2.
\end{eqnarray}
By Cauchy-Schwarz, we obtain
$$ |\eqref{dif1} |  \leq ||\omega - \omega_{MF}||_{L^2} \, ( ||\omega||_{L^2}
 + ||\omega_{MF} ||_{L^2})$$
 and hence we find
 \begin{equation}
 \label{edif1}
  |\eqref{dif1} | \leq  C\,  \| \omega-\omega_{MF} \|_{L^1}^\alpha
  \end{equation}
   for some $\alpha >0$, 
  by using \eqref{pq} with $r=2$ and $q =1$.
 Next, we split \eqref{dif2} as 
 \begin{eqnarray*}
  |\eqref{dif2}| & \leq  & \int (\omega -\omega_{MF}) |\nabla \psi |^2  +
  \int \omega_{MF} \nabla( \psi -\psi_{MF }) \cdot \nabla(\psi + \psi_{MF}) \\ 
   & \leq & ||\omega- \omega_{MF}||_{L^2} \, ||u||_{L^4}^2 \\
    & & \quad  +
     ||\omega_{MF}||_{L^2}\, (||u||_{L^4} + ||u_{MF}||_{L^4} ) ||u - u_{MF}||_{L^4}\, .
  \end{eqnarray*}
  We notice that thanks to \eqref{vomega}, the $L^4$ norm of the velocity
   is bounded in term of the $L^{4 \over 3 }$ norm of the vorticity.
   Therefore, by a new use of \eqref{pq}, we find
  \begin{equation}
  \label{edif2}
   |\eqref{dif2}| \leq  C \, \Big(  ||\omega- \omega_{MF}||_{L^2}
   +  ||\omega- \omega_{MF}||_{L^{4 \over 3}} \Big) \leq C  \| \omega-\omega_{MF} \|_{L^1}^\alpha.
   \end{equation}
  Next, by using   \eqref{edif1}, \eqref{edif2},  we get that
  $$  | D(\om) - D(\om_{MF}) | \leq  C ||\om- \om_{MF} ||_{L^1}^\alpha.$$
   Consequently,  we can use Theorem \ref{theoL1cont} to get that
    \begin{equation}
    \label{Dom}
D(\om ) \geq D(\om _{MF}) - |D(\om ) - D(\om _{MF})| \geq \frac 12 D(\om _{MF})
  \end{equation}
  provided $S(\om ) -  S(\om_{MF} ) $ is sufficiently small.
 Finally, by  using  \eqref{Dom}, \eqref{edif1}, \eqref{edif2} and \eqref{posD},  we easily conclude the proof.  \qed

   \bigskip
   
\section{Global Existence and uniqueness}

We start with a brief  explanation about the  construction of a classical local solution. 
 Let $\om_0 \in L^p$ with $p>2$ be the initial condition such that  $\om_0 \in {\cal M}(E,I)$. Let $\om_{MF}$ the unique Mean-Field solution associated to ${\cal M}(E,I)$
  as above. We assume 
 that $S(\om_0)-S(\om_{MF})$ is   small  enough so that, by virtue of 
 Theorems \ref{theoL1cont}, \ref{theoL1cont of a,b},
   we have
 $$
|D(\om_0) - D(\om_{MF})| \leq Ê{1 \over 2} D(\om_{MF})$$
 and also
 $$ 
|a(\om_0)-a(\om_{MF}) | \leq \frac 12 |a(\om_{MF})|, \,   |b(\om_0)-b(\om_{MF}) | \leq \frac 12 |b(\om_{MF})|.
$$
This implies  that we have the upper bound
$$
|a (\om_0)|+|b(\om_0) | \leq  2  (|a(\om_{MF})| +|b(\om_{MF}) |)
$$
 and also that
$$
 D(\om_{0})  \geq  { 1 \over 2 } D(\om_{MF}) >0.$$
 Note that the  positivity of $D$ is important  in order to  stay away from the singularity.
%  and also that  the non-positivity of $a$ is important  in order to have a confining
%   field.

For every $p\geq 2$,   by using  a standard iterative scheme, 
    we can easily establish, a local
existence and uniqueness result   for a classical solution 
$ \om \in \mathcal{C}([0,T], L^p) \cap L^{\infty}([0, T],  L^1 ( (1+ |x|^2)dx \big)$,   for which
\begin{equation}
\label{upLp}
|| \om (t) ||_{L^p} \leq  2( ||\om_{0}||_{L^p} + ||\om_{MF}||_{L^p}),
\end{equation}
%\begin{equation}
%\lavel{upab}
%|a(\om (t))| \leq 2 |a(\om_{MF}) |, \quad |b(\om (t))| \leq  2|b(\om_{MF}) |, \quad
% \forall t \in [0, T].
%\end{equation}
and
\begin{equation}
\label{lowD}
  \quad D(\om(t)) \geq { 1 \over 4} D(\om_{MF}),
  \quad \forall t\in [0, T].
 \end{equation}
Moreover, we can continue the solution as long as the $L^p$ norm of
 $\om$ remains finite and the denominator $D(\om)$ remains positive.

 Note that $L^2$ seems the natural space for our equation in order
 to have $a$ and $b$ well-defined.
Note also  that the condition \eqref{lowD} allows to avoid the singularity of $D(\om)$
 by Theorem \ref{theoL1cont of a,b}.

Let $T>0$ be the maximal time for which,  the estimates  \eqref{upLp}, 
 \eqref{lowD} are  verified. 
Our purpose is to prove, by a priori estimates,  that $T= + \infty$.
 
 To do this we  shall use the  $L^p$  estimates given by the following Theorem which will be proven in the next section
    \begin{theo}
    \label{theoprop}
    Consider a local  solution as above  of \eqref{nscons} such that
    \begin{equation}
    \label{hypprop}
    |a(\omega(t)) | + |b(\omega(t)) |  \leq C_{0}, \quad \forall t \in [0, T].
    \end{equation}
    Assume also  that   $\omega_{0}\in L^p$, 
    $p\in[2, + \infty)$
     is a probability density.  Then there exists $C_{p}$ which depends only
    on $\omega_{0}$,  $C_{0}$ and $p$   (and hence does  not depend
    on $\nu$ and $T$ if $C_{0}$ does not) such that
   
    \begin{equation}
    \label{estprop}
    ||\omega(t) ||_{L^p} \leq C_{p}, \quad \forall t \in [0, T].
    \end{equation}
    \end{theo}
    
    \bigskip
    
  We  are now  in position to  get a global existence result by showing that $T=+\infty$.  Indeed by  the H-Theorem 
  (see (1.22)), we have  
\begin{equation}
S(\om (t))-S(\om_{MF}) \leq S(\om_0)-S(\om_{MF}), \quad \forall t \in [0, T]
\end{equation}
and  thus  by Theorem \ref{theoL1cont},   we get that $\| \om (t))- \om_{MF}\|_{L^1}
 \leq \eps, \, \forall t \in [0, T]$ 
( with $\eps$ independent of   $T$),  provided that  $S(\om_0)-S(\om_{MF})$ is sufficiently small.  Next,   we can  use   Theorem \ref{theoprop}. Indeed, 
  because of \eqref{upLp}, \eqref{lowD}  we have the bound
  $ |a|+ |b| \leq C_{0}$ 
  on $[0, T]$ for some $C_{0}>0$
  This yields  a control of the $L^p$ norm of $\om$ with $p>2$ , depending only on
   $C_{0}$. In particular, thanks to  \eqref{lowD}, we get
    that
   $$ |D(\om(t)) - D(\om_{MF}) | \leq C \eps^\alpha, \quad
    ||\om(t) - \om_{MF}||_{L^p} \leq C \eps^{\alpha}, \quad \forall t \in [0, T]$$
     where $C$ depends on $C_{0}$ only.
   Consequently, we can choose $\eps$ sufficiently small to have
   $$  ||\om(t) ||_{L^p} < 2 (||\om_{0}||_{L^p}+ ||\om_{MF}||_{L^p}), \quad
    D(\om(t)) > { 1\over 4} D(\om_{MF}), \quad \forall t \in [0, T].$$
and hence  $T=+\infty$.  \qed

Therefore we have proven the following global existence result:
 \begin{theo}
    \label{theo6}
    There exists $\delta_{0} >0$  (independent of $\nu>0$) such that,
     for any initial datum  $\om_0 \in L^p \cap {\cal M}(E,I) $  with $p>2$, 
close to $\om_{MF}$ in the sense that  
    \begin{equation}
    \label{Spetit}
    S(\om_0)-S(\om_{MF}) \leq \delta_{0},
    \end{equation}  
  there exists a unique classical solution
   $$\om (t)
 \in \mathcal{C}([0,+ \infty[, L^p)\cap L^\infty([0, T], \mathcal{M}(E,I))$$ to eq.n \eqref{nscons}  with initial datum $\om_0$.  
 
 Moreover,   we have  the Lyapounov  stability  of $\om_{MF}$, namely,  for any $\eps>0$, there exists
  $\delta$, $0<\delta  \leq \delta_{0}$  such that  if $S(\om_{0}) - S(\om_{MF})
   \leq \delta$, then
 $$ ||\om(t) - \om_{MF}||_{L^1} \leq \eps, \quad \forall  t \geq 0.$$ 
    
\end{theo}  

 As noticed after Theorem
 \ref{theoL1cont}, $\om_{MF}$ is also Lyapounov stable as a stationary solution of the Euler equation
  \ref{euler}. Consequently,  we also have the following global stability result between the flows of
   \eqref{nscons} and \eqref{euler} in the vicinity of $\om_{MF}$: 
 
 \begin{cor}
 For every $\eps>0$, there exists $\delta >0$ such that if $\om_{0} \in L^\infty
  \cap \mathcal{M}(E,I)$ verifies
  $$ S(\om_{0}) - S(\om_{MF}) \leq \delta$$
  then the global solution $\om^E$ of Euler equation \ref{euler}
   and the global solution  $\om^\nu$ of \eqref{nscons} for $\nu>0$ with
    the same initial datum $\om_{0}$ satisfy
    $$ ||\om^\nu (t) - \om^E(t) ||_{L^1} \leq \eps, \quad \forall t \geq 0.$$
 \end{cor} 
 Note that this global approximation property of the Euler evolution  (even for $\nu $ large)
  is of course false for the Navier-Stokes evolution.

\section {  Propagation of $L^p$   regularity }

In this section we prove Theorem \ref{theoprop}.     
            
%I have used a general omega instead of u because we use it also in this sense 

Our stategy in proving Theorem \ref{theoprop} will be based on  weighted energy estimates because
  in this way we can use the fact that the inertial term 
$u \cdot \nabla \omega$ does
   not contribute.  This is crucial in order to find estimates independent of $\nu$.
   
  We first focus on the $L^2$ estimate.  
 We shall use the notation
  $$  \overline{B}=sup_{t\in [0,T]} (|(a(t)|+|b(t)|).$$
  The standard  $L^2$ energy estimate for \eqref{nscons} gives:
   \begin{equation}
  \label{estener1}
  {d \over dt }\Bigl( { 1 \over 2} ||\omega(t)||_{L^2}^2 \Bigr)
 + \nu || \nabla \omega(t)    ||_{L^2}^2 = \nu \,\, \Bigl(  b(t) \int  \omega(t)^3 -  \, 
 a(t) \int \omega(t)^2 \Big).
 \end{equation}
Next, as in \cite{BDP}, \cite{JL}, we can use the Sobolev inequality 
 \eqref{sob2} to get  (we recall that $\int \omega = 1 $)
$$\int \omega^3 = \int (\omega^{ 3 \over 2 })^2 \leq { 9 \over 4}\,  C \Bigr( \int |\nabla \omega | \, 
  \omega^{ 1 \over 2 } \Bigr)^2 \leq { 9 \over 4} C   \,  || \nabla \omega ||_{L^2}^2
  $$
  where $C$ is the best constant in the Gagliardo-Nirenberg-Sobolev inequality.
  Consequently we get
 $$ {d \over dt } \Bigl({ 1 \over 2} ||\omega(t)||_{L^2}^2 \Bigr)
 + \nu || \nabla \omega(t)    ||_{L^2}^2  
\leq \nu  C\,  \overline{B} \,  ||\nabla \omega(t) ||_{L^2}^2 
  + \nu\,  \overline{B}\,  ||\omega(t)||_{L^2}^2
   $$
   where $C$ is an explicit harmless number.
Now,  let us assume for the moment that  $C  \overline{B} $
     is sufficiently small (less than $1/2$ for example),  then we can deduce that
 \begin{equation} 
 \label{ener1}
  {d \over dt } \Bigl({ 1 \over 2} ||\omega(t)||_{L^2}^2 \Bigr)
 + {\nu \over 2 }  || \nabla \omega(t)    ||_{L^2}^2  \leq \nu 
 \,   \overline{B}\,  ||\omega(t)||_{L^2}^2.
 \end{equation}
If we directly integrate this differential  inequality, 
we still  cannot conclude. Indeed we shall
 find that $||\omega(t)||_{L^2}$ grows exponentially in time and this
 does not allow us to get a  uniform in time  estimate.  Note however that the bad term in the right hand-side of (4.7) comes from the
 linear   term $ \nabla\cdot ( x \, \omega)$ in  eq.n
 \eqref{nscons}. The explanation for the bad behaviour we get  is simple: 
  the semigroup generated
 by the linear  Fokker-Planck operator
     $$ L \omega = \Delta \omega + \overline{B} \nabla \cdot (x \omega),$$
     with $\overline{B}>0,$ is not uniformly   bounded   in time 
as an operator in $\mathcal{L}(L^2).$
  Nevertheless, it is bounded as an operator 
in $\mathcal{L}(L^1 \cap L^2, L^2)$.
 A very simple way to see this property is to use a weighted energy estimate.
  Indeed, multiplying \eqref{estener1} by 
$e^{\nu\, \overline{B}\, t}$, we find
\begin{equation}
\label{ener2}
 {d \over dt }\Bigl( { e^{\nu \, \overline{B}\,  t} \over 2} ||\omega(t)||_{L^2}^2 \Bigr)
 + {\nu \, e^{\nu\, \overline{B} t} \over 2 }|| \nabla \omega(t)    ||_{L^2}^2 \leq
   { 3 \over 2 }\, \nu \, \overline{B}\,
  e^{\nu \overline{B}  t } ||\omega(t) ||_{L^2}^2
 , \quad \forall t \in [0, T].
 \end{equation}
Now,  we can use  \eqref{sob2} and the Young inequality to get,
 for some harmless explicit number $C$ (independent of $\nu$) which
  changes from line to line,
$$ { 3 \over 2 }\, \nu \, \overline{B}\,
  e^{\nu \overline{B}  t } ||\omega(t) ||_{L^2}^2 
  \leq  C\, \nu \, \overline{B}\, e^{\nu \overline{B}  t }
   ||\nabla \omega(t) ||_{L^2}  \leq { 1 \over 4 } \, \nu \,e^{\nu \overline{B}  t }
    ||\nabla \omega(t) ||_{L^2}^2  + C \, \nu  \, \overline{B}^2\,  \,e^{\nu \overline{B}  t }. 
    $$
    Consequently, we can plug this last inequality in \eqref{ener2}
     to get 
 $$ {d \over dt }\Bigl( { e^{\nu \, \overline{B}\,  t} \over 2} ||\omega(t)||_{L^2}^2 \Bigr)
   \leq C \, \overline{B}^2\,  \nu  \,e^{\nu \overline{B} \, t }. $$
The integration finally gives
$$ ||\omega(t)||_{L^2}^2 \leq ||\omega_{0}||_{L^2}^2  + C\,  \overline{B},
 \, \quad \forall t \in [0, T ].$$
 
 We now remove the assumption on the smallness of  $\overline{B}$ and prove also the propagation of the $L^p$ regularity in the general case  $p\geq2$.

  The starting point is to use  the idea of \cite{JL}, \cite{BDP}
   in the study of the Keller-Segel equation. For $K > 1,$  a parameter
which will be fixed later, we define
   $$ m_{K}(t) = \int (\omega(t) -K)_{+}\, dx\,.$$
We note that
   \begin{equation}
   \label{mK0}
   m_{K}(t)   \leq
    \int_{\omega(t)\geq K} \omega (t)  \leq { 1 \over \log K } \int \omega(t) \log(\omega(t))
     \, dx  
    \leq { 1 \over \log K} \int \omega(t) |\log \omega(t) |.
   \end{equation}
   Now we can use the following useful   inequality : 
   
   \begin{lem}
   \label{lement}
 There exists $C>0$,   such that,  for all   probability distribution $\omega$, we have
  \begin{eqnarray}
  \label{entropie2}
  \int \omega |\log \omega | \leq
   S(\omega)  + C(1 +  I(\omega))
    \end{eqnarray}
 
\end{lem}

We postpone the proof of the lemma to the end of the section.

\bigskip

 Thanks to \eqref{entropie2}, we
find that  $\int \omega |\log \omega |$ is bounded in terms of the 
initial datum because the entropy is decreasing and $I(\om)$ is constant. Thus
\begin{equation}
\label{mK}
 m_{K}(t) \leq { 1 \over \log K } C(\omega_{0}), \quad \forall t \geq 0.
 \end{equation}
  
      Next, we  can perform a modified $L^p$ energy estimate for the solution
       of  \eqref{nscons}. After
       a few integrations by parts, we find
 \begin{eqnarray}
 \label{p1}
 & &   
    { d \over dt }\Big(
     {1 \over p } \int (\omega-K)_{+}^p \Big)
 + \nu (p-1) \int  (  \omega - K )_{+}^{p-2}  |\nabla (\omega-K)_{+}  |^2 \\
\nonumber 
& &   = \nu \overline {B}  \Big( \int  (\omega - K )_{+}^{p+1} + ( 2K + { 2 \over p } ) + ( K^2 + 2 K)
   \int ( \omega - K)_{+}^{p-1} \Big).
   \end{eqnarray}
   To estimate the first term in the right hand side of \eqref{p1}, we use 
     the Sobolev-Gagliardo-Niremberg inequality \eqref{sob1}
      and Cauchy-Schwarz. We have
  \begin{eqnarray*}
   \int  (\omega - K )_{+}^{p+1} = \int \Big(  (\omega - K )_{+}^{p+1 \over 2} \Big)^2
  &  \leq &  C\, {(p+1)^2 \over 4}\,   \Big( \int (\omega-K)_{+}^{p-1 \over 2 } |\nabla ( \omega - K )_{+} | 
   \Big)^2  \\
 &   \leq  & C\,  {(p+1)^2 \over 4}\, m_{K}(t)  \int  (\omega- K )_{+}^{p-2}
    |\nabla ( \omega - K )_{+} |^2.
\end{eqnarray*}
This yields,  thanks to our assumption \eqref{hypprop} and \eqref{mK},
$$ \nu \, \overline {B} \,  \int  (\omega - K )_{+}^{p+1} \leq\nu\,  {C_{0} \,S_{0}\, (p+1)^2 \over 4 \log K } 
  \int  (\omega- K )_{+}^{p-2}
    |\nabla ( \omega - K )_{+} |^2.$$
    Hence, by choosing  $K$  such that
    $$ {C_{0} \,S_{0}\, (p+1)^2 \over 4 \log K } = { 1 \over 2 } (p-1), $$
  we obtain  
\begin{equation}
\label{p+1est}
 \nu \,  \overline {B} \,  \int  (\omega - K )_{+}^{p+1} \leq { \nu\, (p-1) \over 2 } 
  \int  (\omega- K )_{+}^{p-2}
    |\nabla ( \omega - K )_{+} |^2.
 \end{equation} 
 Note that $K$ depends only on $C_{0}$,  $S_{0},$  $p$ and is 
diverging with $p$.

 To estimate the last term in \eqref{p1}, we write
 \begin{eqnarray}
 \nonumber
  \int (\omega - K )_{+}^{ p - 1 }&  \leq & 
  \int_{K \leq \omega \leq K+1}( \omega- K )_{+}^{p-1} + 
   \int_{ \omega \geq K+1 } ( \omega - K )_{+}^{p-1} \\
 \label{p-1est}
 & \leq &  1 + \int (\omega - K )_{+}^p
  \end{eqnarray}
 By plugging \eqref{p-1est} and \eqref{p+1est} in \eqref{p1}, we find
 \begin{eqnarray*}
   & &{ d \over dt }\Big(
     {1 \over p } \int (\omega-K)_{+}^p \Big)
 +{ \nu (p-1) \over 2 } \int  (  \omega - K )_{+}^{p-2}  |\nabla (\omega-K)_{+}  |^2 \\
\nonumber
& & \leq \nu C_{0}\Big( K^2  + 4 K +  { 2 \over p }  \int (\omega - K)_{+}^p + K^2+
 2 K\Big) \\
\nonumber 
 & & \leq C_{0}\, C\, \nu \Big( \int (\omega - K)_{+}^p +1 \Big)
 \end{eqnarray*}
 where,  from now on,  $C$ is a harmless number which depends only on
  $K$ and $p$.  Again, we note that we cannot directly conclude by
   using the Gronwall Lemma in the last differential inequality because 
   it  gives an estimate which is not uniform in time. We  now 
   use the technique that we have explained in the beginning.
  We find
\begin{eqnarray}
\label{p2}
   & &{ d \over dt }\Big( e^{\nu \, t }
     {1 \over p } \int (\omega-K)_{+}^p \Big)
 + {\nu (p-1) \over 2 }
  e^{\nu \,t } \int  (  \omega - K )_{+}^{p-2}  |\nabla (\omega-K)_{+}  |^2 \\
 \nonumber
 & & \leq C ( 1 +  C_{0}) \, \nu \,  e^{ \nu\, t }\, \Big( \int (\omega - K)_{+}^p +1 \Big).
 \end{eqnarray}
 Next, we can use the inequality \eqref{sob2} to get
 \begin{equation}
 \label{sobp1} \int (\omega - K)_{+}^p = 
  \int\Big(  (\omega - K)_{+}^{p \over 2 } \Big)^2 \leq C \int ( \omega - K )_{+}^{p \over 2}
  \,\Big(  \int ( \omega - K)_{+}^{p-2  } |\nabla (\omega - K )_{+} |^2 \Big)^{ 1 \over 2}.
 \end{equation}
 By  using the interpolation inequality \eqref{pq} of $L^{p/2}$ between
 $L^1$ and  $L^p$, we  have 
$$ \int  ( \omega - K)_{+}^{p \over 2  } \leq C \Big( \int ( \omega -K )_{+}^p \Big)^{p-2
 \over 2(  p-1) } $$
 and hence, we deduce from \eqref{sobp1} that
 $$  \int ( \omega - K )_{+}^p \leq C \Big( \int (\omega - K)_{+}^{p-2}
  |\nabla(\omega - K )_{+}|^2 \, \Big)^{ 1 \over q}$$
  where $q$ is such that $p^{-1} + q^{-1 } = 1$. Thanks to the 
   inequality 
   \begin{equation}
   \label{cvx}
    ab \leq p^{-1} a^p + q^{-1 } a^q,  \quad a \geq 0, \, b \geq 0, 
    \end{equation}
   we finally obtain
  $$ C( 1 + C_{0}) \int ( \omega - K )_{+}^p \leq { p-1 \over 4 }
  \int (\omega - K)_{+}^{p-2}  |\nabla(\omega - K )_{+}|^2 + C_{p},$$
  where $C_{p}$ will now  stand for a number which depends only
   on $\omega_{0}$, $C_{0}$
   and $p$. 
  By using this last inequality in \eqref{p2}, we finally arrive to 
  $$ { d \over dt }\Big( e^{\nu \, t }
     {1 \over p } \int (\omega-K)_{+}^p \Big) \leq  C_{p}\, \nu e^{\nu t }.
   $$
   The integration gives
   \begin{equation}
   \label{Kp+}
    \int (\omega -K)_{+}^p \leq C_{p}.
   \end{equation}
   Now we can conclude as
    in \cite{JL}, \cite{BDP}.  By using the inequality
     $$ x^p \leq \Big( {\lambda \over \lambda -1} \Big)^{p-1} (x-1)^p$$
      for every $x \geq \lambda >1$, we find
    \begin{eqnarray*}
     \int \omega^p & \leq &
          \int_{\omega \leq K} \omega^p + 
           \int_{\omega > K } \omega ^p \\
      &\leq  & K^{p-1} + \int_{K<\omega \leq \lambda K } \omega^p
              + \int_{ \omega \geq \lambda K} \omega^p \\
      & \leq & K^{p-1} + (\lambda K)^{p-1}  +
        K^p\Big( {\lambda \over \lambda -1} \Big)^{p-1}
              \int_{\omega \geq \lambda K } \Big( {\omega \over k } - 1 \Big)^p  \\
      &\leq  & K^{p-1} + (\lambda K)^{p-1} + \Big( {\lambda \over \lambda -1} \Big)^{p-1}
              \int (\omega - K ) _{+}^p.
    \end{eqnarray*}
   This ends the proof of Theorem \ref{theoprop}   \qed.

   \bigskip
   
   It remains to prove Lemma \ref{lement} which is a classical estimate we present for completeness.
   
   Define $\overline{\omega} = \omega {\bf 1}_{\{|\omega 
 | \leq 1 \} } $ Since we have
\begin{equation}
\label{lem1est1} \int \omega |\log \omega | =  S(\omega)  - 2 \int \overline{\omega} \log \overline
{\omega}, 
 \end{equation}
 it suffices to find a bound from below of 
$\int \overline{\omega} \log \overline
{\omega}$.
  By using the fact that the relative entropy between two probability measures is non negative (this is an easy consequence of the Jenssen inequality) we get 
 $$ \int_{\mathbb{R}^2} (\overline{\om}/{\overline{m}})\log\Big( {
  \overline{\om}/{\overline m} \over  { 1 \over 2 \pi} e^{- {|x|^2 \over 2 } } }\Big)\geq 0$$
  where $\overline{m}=\int_{\mathbb{R}^2} \overline{\om}\leq 1.$
%$$
% \int \overline{\omega} \Bigl(\log \overline{\omega} + { 1 \over 2} |x|^2\Bigr)
%  =   \int \overline{\omega} \log \Bigl( {  \overline{\omega} \over  { 1 \over 2 \pi } e^{ -   |x|^2
%  \over 2 }  } \Bigr)+ \int \overline{\omega} \log \left( { 1 \over 2 \pi } \right) 
%   = \int \overline{\omega} \Bigl(\log \Bigl( {  \overline{\omega} \over  { 1  \over 2 \pi } e^{ -   |x|^2
%  \over 2 }  } \Bigr) +\log \left( { 1\over 2 \pi } \right)\Bigr)
%  $$
%  Let us set $\Omega = \overline{\omega}/  { 1  \over 2 \pi } e^{ -   |x|^2
%  \over 2 }$, $d\mu = { 1 \over 2 \pi } e^{ -  |x|^2
%  \over 2 } dx, $ then by Jenssen inequality
% $$ \int \overline{\omega} \log \Bigl( {  \overline{\omega} \over  { 1 \over 2 \pi } e^{ - |x|^2
%  \over 2 }  } \Bigr) = \int \Omega \log \Omega \,d\mu \geq \Bigl( \int \Omega \,d\mu \Bigr)
%   \log \Bigl( \int \Omega \, d\mu\Bigr) $$
%    since 
%    $ \int \Omega\,  d\mu = \int \overline{\omega}\, dx \leq 1,$ we have
%    $$ \Bigl( \int \Omega \,d\mu \Bigr)
%   \log \Bigl( \int \Omega \, d\mu\Bigr) \geq  - e^{-1 }.$$
   Then  we get 
    $$ 
    \int \overline{\omega} \log \overline{\omega} \geq
    -  \overline{m}I(\omega) + \overline{m}\log{ 1 \over 2\pi} -\overline{m}\log\overline{m}\geq 
    -  I(\omega) + \log{ 1 \over 2\pi} -\frac{1}{e}
    $$
    and hence we get \eqref{entropie2} by using this last estimate and \eqref{lem1est1} \qed

   \section{Asymptotic behaviour  }
    \label{asymptotic}
    
   In this section we investigate the asymptotic behaviour of the global solutions
    given by Theorem \ref{theo6}.  More generally, one can consider
     a global solution $\om(t)$ of \eqref{nscons} such that
      $\om\in \mathcal{C}([0, + \infty[, L^2 \cap L^1( (1+ |x|^2) dx)$
       and such that $\om(t) \in \mathcal{M}(E,I)$ and which satisfies
           for some $C >0,$     the uniform estimates
    \begin{equation}
      \label{ashyp}
       ||\omega (t ) ||_{L^2} \leq C, \quad  |b(\omega(t))|+|a(\om(t))|
      \leq C, \quad \forall t \geq 0.     
     \end{equation}
    
    The main result of this section is given by the following Theorem.
    
      \begin{theo}
  \label{theoas}
  Let $\omega(t)$ a  global solution  of \eqref{nscons} as above 
   which satisfies  \eqref{ashyp}.   
   Then $\omega(t)$ converges in $L^1$, as $t\to \infty$, to the unique
   solution $\omega_{MF} \in { \cal M}(E,I)$   of  the associated microcanonical variational problem.  
  \end{theo}
  
  Note that the solutions constructed in Theorem \ref{theo6} satisfy the estimate
   \eqref{ashyp} and hence  their asymptotic behaviour is given by Theorem
    \ref{theoas}.

  \subsection*{Proof of Theorem \ref{theoas}}
The first step consists in  proving  that the orbit $\{\omega(t)\}_{t \geq 0 }$  is
 relatively   compact
 in $L^1$ and uniformly  bounded in $L^\infty$.  Before,  we  need  to study the 
 evolution operator  generated
  by the non-autonomous Fokker-Planck type operator 
  $$ L_{\gamma} \omega =  \nu \Bigl( \Delta \omega + \gamma(t) \nabla\cdot ( x \omega) 
   \Big)$$
  where $\gamma(t)$ is a given continuous curve. Denoting by $S_{\gamma}(t, \tau)
  \om_{0}$
   the solution of
   $$ \partial_{t} \om = L_{\gamma} \om,\quad t>\tau,  \quad \om(\tau)= \om_{0},$$
     we have the following estimates
  \begin{lem}
  \label{lemFP}
  Suppose that for all $t \geq 0$, $ |  \gamma(t) |  \leq K_{0}$, for some $K_{0}>0$.
  Then, there exists $C>0$ independent of  $\nu >0$ and   such that 
   for $(p,q,r) \in [1, + \infty]^3$  we have: 
   \begin{eqnarray}
   \label{FP1}
  &&   || S_{\gamma}(t, \tau) \omega ||_{L^p} \leq C  K_{0}^{ 1 - { 1 \over r } } 
  { e^{  2 \nu K_{0}( 1 - { 1 \over q } ) (t- \tau) }
     \over ( 1 - e^{ - 2 \nu  K_{0}( t - \tau) })^{ 1 - { 1 \over r } } }
      ||\omega ||_{L^q}, \quad { 1 \over r } + { 1 \over q} =  1 + { 1 \over p },  \\
 \label{FP2}
  && ||\nabla S_{\gamma}(t, \tau) \omega  ||_{L^p}  
     \leq   
 C  K_{0}^{ { 3 \over 2 } - { 1 \over r } } 
  { e^{  2 \nu K_{0}( 1 - { 1 \over q } ) (t- \tau) }
     \over ( 1 - e^{ - 2 \nu  K_{0}( t - \tau) })^{  {3 \over 2 } - { 1 \over r } } }
      ||\omega ||_{L^q}, \quad { 1 \over r } + { 1 \over q} =  1 + { 1 \over p } \\
\label{FP3}
  && || S_{\gamma}(t, \tau) \nabla \omega  ||_{L^p}  
     \leq   
 C  K_{0}^{ { 3 \over 2 } - { 1 \over r } } 
  { e^{   \nu K_{0}( 1 - { 2 \over q } ) (t- \tau) }
     \over ( 1 - e^{ - 2 \nu  K_{0}( t - \tau) })^{  {3 \over 2 } - { 1 \over r } } }
      ||\omega ||_{L^q}, \quad { 1 \over r } + { 1 \over q} =  1 + { 1 \over p }        
  \end{eqnarray}  
  for all $p\in [1, + \infty]$
  \end{lem}
  
  \subsubsection*{Proof of Lemma \ref{lemFP}}
  A simple computation in Fourier space allows  to find the 
explicit representation
\begin{equation}
S_{\gamma}(t, \tau) \omega (x)=
   e^{  2 \nu ( B(t) - B(\tau)  ) }
\int
\frac
{  e^{ - \tfrac{ | x - y |^2}{4 \nu \int_{\tau}^t  e^{ 2 \nu  (B(s) - B(t)) }\,ds  }}}
{4 \pi \nu \int_{\tau}^t  e^{  2 \nu (B(s) - B(t)) } \, ds }    
\omega_{0}( e^{ \nu( B(t) - B(\tau)  )}y ) \, dy
\end{equation}
where $B(t)= \int_{0}^t \gamma (s) \, ds.$ The result of Lemma \ref{lemFP}
then follows by standard convolution estimates. \qed

\bigskip

We come back to the proof of Theorem \ref{theoas}. We shall prove that
 $ ||\nabla \omega(t) ||_{L^1}$ is uniformly bounded. We use the same idea
 as in \cite{GW}.
  Note that the solution of  \eqref{nscons} can be written as
  \begin{equation}
  \label{asduhamel}
  \omega(t) = S_{a}(t,0) \omega_{0}
   + \int_{0}^t S_{a}(t,\tau)\nabla \cdot  \Big( - u \,\omega + \nu b u^\perp\,
    \omega
    \Big)(\tau) 
    \, d\tau.
   \end{equation}
 Moreover, thanks to \eqref{ashyp},   we  have a uniform estimate on $ ||\omega (t)||_{L^2}$
     and on $|a(t)|+|b(t)|$ for all times : there exists $C_{0}>0$ such that
     \begin{equation}
     \label{asbound}
     ||\omega(t) ||_{L^2} \leq C_{0},  \quad  |a(t)|+|b(t)| \leq K_{0}, \, \forall t \geq 0.
     \end{equation}
Consequently, in \eqref{asduhamel}, we can consider
 $a$ and $b$ as known and  we can  use the estimates of Lemma \ref{lemFP}.
 Let us define
 $F(\omega)$
   as the right-hand side of \eqref{asduhamel}. We have 
       \begin{eqnarray*}
     ||F(\omega(t))||_{L^\infty}
     & \leq  &  C \Big( {C_{0} e^{ \nu K_{0} t  }   \over a_{\nu}(t)^{ 1 \over 2  } } 
      +  C_{0}(1+ \nu ) \int_{0}^t   {e^{ \nu K_{0}( t- s)  } \over a_{\nu}(t-s)^{ 3 \over 4} }
       ||u \,  \omega ||_{L^4} \, ds \Bigr), 
   \end{eqnarray*}
   where $a_{\nu}(t) = 1 - e^{ - 2 \nu K_{0} t }. $
   Since by \eqref{vomega},  \eqref{pq} and the uniform $L^2$ bound, we have  : 
   $$ 
   ||u \, \omega ||_{L^4} \leq   ||\omega ||_{L^\infty} \, ||u ||_{L^4}
    \leq C  ||\omega ||_{L^\infty}  \, ||\omega ||_{L^{ 4 \over 3 } } \leq C C_{0} ||\omega ||
    _{L^\infty}, 
    $$
    we finally get
  $$  ||F(\omega(t))||_{L^\infty}
      \leq    C \Big( {C_{0} e^{ \nu K_{0} t  }   \over a_{\nu}(t)^{ 1 \over 2  } } 
      +  C_{0}(1+ \nu ) \int_{0}^t   {e^{ \nu K_{0}( t- s)  } \over a_{\nu}(t-s)^{ 3 \over 4} }
       ||\omega(s)||_{L^\infty} \, ds \Big).$$
       Consequently,  we can set
       $$ z(T) = \sup_{[0, T ]} \Big(e^{ - \nu K_{0} t } a_{\nu}(t)^{ 1 \over 2} ||\omega(t)||_{L^
       \infty}\Big)$$
       to get 
       $$ z(T) \leq C C_{0}(  1 + \nu) \Big( 1 + a_{\nu}(T)^{ 1 \over 2} 
        \Big( \int_{0}^T { 1  \over a_{\nu}(T-s)^{ 3 \over 4} \, a_{\nu}(s)^{ 1 \over 2}  } \, ds \Big) z(T)\Big).
       $$
       Next, we notice that 
       $$ \lim_{T \rightarrow 0}a_{\nu}(T)^{ 1 \over 2}
        \int_{0}^T  { 1  \over a_{\nu}(T-s)^{ 3 \over 4} \, a_{\nu}(s)^{ 1 \over 2}  } =0,$$
therefore, there exists $T(\nu, C_{0})>0$ such that 
$$ z(T(\nu, C_{0})) \leq CC_{0}( 1+ \nu ) + { 1 \over 2 } z(T(\nu, C_{0})) $$
 and hence, we get that
 \begin{equation}
 \label{aslinfty1}
  ||\omega (t) ||_{L^\infty} \leq { C(\nu, C_{0})  \over a_{\nu}(t)^{ 1 \over 2} }, \quad 
 \forall t \in [0, T(\nu, C_{0})]\,.
 \end{equation} 
Next, since to establish  \eqref{aslinfty1} we have only used \eqref{asbound},
 we can consider for every $n \in \mathbb{N}$, 
  the solution $\tilde{\omega}$ of \eqref{nscons} with initial value
  $\omega(n T(\nu, C_{0})/2)$. By the above argument, we
   get that $\tilde{\omega}$ satisfies the estimate \eqref{aslinfty}.
    By uniqueness, we  have
  $$ \tilde{\omega}(t) = \om( t + n T(\nu, C_{0})/2 ) , \quad \forall t \in [0, T(\nu, C_{0})]$$
 and hence, 
    \begin{equation}
  \label{aslinfty}
  ||\omega (t) ||_{L^\infty} \leq  C(\nu, C_{0})\Big( 1 + { 1   \over a_{\nu}(t)^{ 1 \over 2} }
  \Big), \,
 \forall t \geq 0. 
 \end{equation}
 for some $C(\nu, C_{0})$.
 \bigskip

   In a similar way  we have  by Duhamel formula
\begin{equation}
\label{asduhamelgrad}
 \nabla \omega(t) = \nabla S(t,0) +  \omega_{0} \int_{0}^t \nabla S(t-\tau)
 \Big( - u \cdot \nabla \omega  +   \nu b  u^\perp \cdot \nabla \omega +
 \nu  b  \omega^2 \Big)(s) \, ds
  \end{equation}
   $$ ||\nabla F(\omega) ||_{L^1}
    \leq C \Bigl(  { C_{0} \over  a_{\nu}(t)^{ 1 \over 2 } } + (1+ \nu) C_{0}
     \int_{0}^t { 1 \over a_{\nu}(t-s )^{ 1 \over 2 }  }
      (||u  \,  \nabla \omega (s)  ||_{L^1} + C_{0}) \, ds \Big)$$
    and since we have by \eqref{vinfty} 
    $$ ||u \,\nabla  \omega ||_{L^1} \leq C || u ||_{L^\infty} ||\nabla \omega ||_{L^1}
    \leq C ||\omega ||_{L^\infty}^{ 1 \over 2 }\, ||\nabla \omega ||_{L^1},$$
  we get, thanks to \eqref{aslinfty},
  $$ ||\nabla F(\omega) ||_{L^1}
    \leq C \Bigl(  { C_{0} \over  a_{\nu}(t)^{ 1 \over 2 } } + C_{0}\int_{0}^t
     { 1 \over a_{\nu} (t- s)^{ 1 \over 2}  }\, ds + C(\nu, C_{0})
     \int_{0}^t { 1 \over a_{\nu}(t-s )^{ 1 \over 2 }   a_{\nu}(\tau)^{ 1 \over 4} } \, ||\nabla
     \omega(s)||_{L^1}\, ds \Big).$$
     Consequently, by using the same method as before, 
       we   can easily  obtain    
  \begin{equation}
  \label{asnablaomegaL1}
  ||\nabla \omega(t) ||_{L^1} \leq  C(\nu, C_{0}) \Big( 1 + { 1  \over a_{\nu}(t)^{1 \over 2 } }\Big), \, 
  \forall t \geq 0
  \end{equation}
  for some $C(\nu, C_{0})$.

  \bigskip

  We now consider  $\Omega$,  the omega limit set  of the  trajectory $(\omega(t))_{t \geq 0}.$ We deduce  from the previous estimates that the positive orbit $\{\omega(t)\}_{
  t \geq 0 }$
    is relatively compact in $ X = L^1 ( ( 1 +  |x|^\alpha)\, dx) 
    \cap L^2$ for $\alpha <2$. Indeed,
    since $\omega(t) \in \mathcal{C}([0, + \infty[, X)$, it suffices to
     prove that   $\{ \omega(t) \}_{t\geq 1}$ is relatively compact.
     The compactness in $L^1( 1+ |x|^\alpha)$ follows 
immediately from the Riesz-Frechet-Kolmogorov criterion: 
      $\omega(t)$ is uniformly bounded in $L^1$, 
      the uniform (for $t\geq 1$)  bound \eqref{asnablaomegaL1}  gives
     the equi-integrability and  we have a uniform bound on the moment 
of inertia
     to control the mass far away. Next,  thanks to the uniform $L^{\infty}$
     estimate for $t \geq 1$ given by \eqref{aslinfty} and the relative
      compactness in $L^1$, we also get that $\{\omega(t)\}_{t \geq 0 }$ is 
       relatively compact in $L^p$ for every $p<+ \infty$.

  By the relative compactness properties that we have just proven, we
  get that $\Omega$ is non empty and actually made by smooth $L^p$
   functions thanks to the smoothing effect of the parabolicity. 
Moreover, 
    we also have  that the elements of $\Omega$
    are probability densities.
    %such that $M(\omega) = \int x \omega = 0.$ 
Also,  if $\omega \in \Omega$,
    since there exists an increasing sequence $t_{n}$ such that
    $\omega(t_{n})$ tends to $\omega$ in $X$, we also have
      \begin{equation}
  \label{asei}
  E(\omega)= \lim_{n}E(\omega(t_{n}))=E, \quad I(\omega) 
\leq  \lim_{n} I(\omega(t_{n}))=I.
  \end{equation}  
    The first equality is proven in the Appendix, see \eqref{En}. Finally, we notice that 
    the entropy  $S$ is constant on $\Omega$. 
Indeed, if $\omega_{1}, \, \omega_{2}\in \Omega$, 
      we can construct  an increasing sequence $t_{n}$ such that
     $\omega(t_{2n}) $ tends to $\omega_{1}$ and $\omega(t_{2n+1})$
      tends to $\omega_{2}$  almost everywhere and such that  there exists
      $g_{1}, \, g_{2} \in 
            L^1((1 + |x|^2)dx) \cap L^2$ with
            $$ \omega(t_{2n}) \leq g_{1}, \quad \omega(t_{2n+1})
             \leq g_{2}.$$
     By using that 
     $$ \omega |\log \omega | \leq C\Big(  \om^2 +  |\om|^{ 3 \over 4}
      \Big) \leq C   \Big(  \om^2 +  ( 1 + |x|) \om + { 1 \over ( 1 + |x|  )^{3} } \Big), $$ 
      we find by Lebesgue Theorem that 
      $$ S(\omega_{1}) = \lim_{n}S(\omega(t_{2n})), \quad
       S(\omega_{2}) = \lim_{n}S( \omega( t_{2n+1} ) ).$$ 
        But,  since the entropy is decreasing, we  also have
       $S(\omega(t_{2n})) \geq S(\omega(t_{2n+1})) \geq S(\omega(t_{2n+2}))$
        so that 
       passing to the limit,   we  get
     $$ S(\omega_{1}) \geq S(\omega_{2}) \geq S(\omega_{1})$$
     and hence $S(\omega_{1}) = S(\omega_{2}).$
    
    Finally, we can prove that the elements of $\Omega$ are solutions 
of the mean field equation. If $\omega \in \Omega$,  consider
      $\omega(t)$ the solution of \eqref{nscons} with initial value $\omega$.  
By the strong parabolic principle, we have that $\omega(t)$ is smooth 
and strictly positive
       for $t>0$. Since $\Omega$ is invariant and $S$ is constant on it, 
the entropy dissipation   identity (1.22) gives that for $t>0$
  $$ \nabla\Big( \log \omega(t) - b \psi(t)-  a  { |x |^2 \over 2 } )  \Big) = 0.$$
          By continuity in time we  get that
          $\omega$ actually solves the mean field equation
         \begin{equation}
         \label{asMF}
         \omega = { 1 \over Z} e^{ b  \psi  + a  { |x|^2 \over 2 } }
         \end{equation}
          in $\mathbb{R}^2.$ Finally,  by using the result of \cite{Naito}
            and
            Lemma 4.3 of \cite{BDP}, we get that $\omega$ is radially symmetric.
             Note that, we also necessarily  have that  $a<0$ and $b<8 \pi.$
       \bigskip
           
        To summarize,  we have proven  that the omega limit set
$\Omega$ of  $\{\omega(t)\}_{t\geq 0}$ is made by probability densities
 which are radially symmetric  solutions of the
          mean field equation \eqref{asMF} with finite energy equals to $E$  and finite  moment of
          inertia. Since the entropy separates the radial mean field solutions
          (see Remark \ref{remap} in the Appendix), 
           we conclude that $\Omega$ consists  in a single point. \qed

 \appendix

\section{The Mean-Field Equation and related variational problems}

In this Appendix we collect  some useful facts concerning the 
Mean-Field Equation (MFE) in $\mathbb{R}^2$. 
The main ideas are in \cite{CLMP1}, \cite{CLMP2}, 
here, we adapt the results to the $\mathbb{R}^2$ case.  
For the microcanonical problem,  the strategy of 
the proof  is slightly different,  we do not  prove directly 
the existence of  a solution. We focus on the negative
 temperature case (which corresponds to $b>0$) which is the most interesting
  case.  
\begin{defi}[Canonical Variational Principle]

For $a<0$, $b>0$, consider the
free-energy functional
\begin{equation}
\label{freeab}
F_{a,b}(\om)=S(\om)-bE(\om)-aI(\om)
\end{equation}
defined on the space  $\Gamma$ of   probability densities on $\Bbb R^2$ for which
$E(\omega),$ $I(\omega)$ and $S(\omega)$ are finite. 
  We set 
 \begin{equation}
\label{free}
F(a,b)=\inf_{\om\in\Gamma}F_{a,b}(\omega)
\end{equation}
\end{defi}
\begin{defi}[Microcanonical Variational Principle]

 For  $E \in \mathbb{R}$ and $ I >0$ let us define 
 $$\mathcal{M}(E,I)=\left\{\omega\in\Gamma:E(\om)=E,I(\om)=I\right\}.$$
 We set 
  \begin{equation}
\label{entroprin}
S=\inf_{\om\in\mathcal{M}(E,I)}S(\omega)
\end{equation}
\end{defi}

The main results of this Appendix are the  two following theorems.
For the positive temperature case (i.e. $b<0$),   $F_{b,a}$ is a convex
 functional so  that there  is    a unique minimizer. Moreover, there exists a unique 
  solution to the MFE  \cite{Go-Lions} and  all the following results are obvious.

 \begin{theo}[Canonical Variational Principle]
 \label{theoMFE}
 For $a<0,$ and $ 0< b<8\pi$ : 
 \begin{itemize}
 \item[i)] there exists $\omega \in \Gamma$ such that 
 $F_{a,b}(\omega) =  F(a,b) $. 
 Moreover $\omega$ is radially symmetric and
  solves the mean field equation \eqref{MFeq}.
  
  \item[ii)]  There is only one radially symmetric solution of the mean field equation
   \eqref{MFeq}.
  
  \item[iii)] As a consequence,  there exists a unique minimizer 
$\omega_{a,b}$ of $F(a,b)$ over $\Gamma$.
  \end{itemize}
  
  \end{theo}
  
  \begin{rem}It is easy to prove that when $b\to 8\pi$, the solutions to the MFE
concentrates at the origin. Indeed by multiplying the equation by $x\cdot
\nabla \psi$ and integrating by parts, we arrive to the identity
(that is the same argument leading to the Pohozhaev inequality): 
\begin{equation}
1- \frac {2\pi a}{b} I=\frac {8\pi}{b}.
\end{equation}
Hence when $b\to 8\pi$, $I\to 0$ and the concentration takes place. For
$b>8\pi$ we do not have solutions.
\end{rem}

As a consequence,  we can solve  the microcanonical variational principle.  
 \begin{theo}[Microcanonical Variational Principle]
 For $a<0,$ and $ 0< b<8\pi$ : 
 \label{MFEmicro}
\begin{itemize} 
\item[i)] F(a, b) is a concave smooth function. 

\item[ii)]  \begin{equation}\frac{\partial F}{\partial a}= - I(\omega_{a,b})<0,
\frac{\partial F}{\partial b}=-E(\omega_{a,b})
\end{equation}

\item[iii)]  For $I>0,E \in \mathbb{R}$ let us define $S^*(I,E)$ as
\begin{equation} 
S^*(I,E)=\sup_{a,b}\Big( F(a,b) + b E + a I \Big) ,
\label{Lgdr}
\end{equation}

Denote by  $a(I,E),b(I,E)$   the unique maximizer  of \eqref{Lgdr},
then for any $I>0, E\in \mathbb{R}$, $S(I,E)=S^*(I,E),$ and hence $S$  is a smooth convex function. Moreover, 
 the microcanonical variational 
principle admits a unique minimizer $\tilde\omega_{I,E}$ in $\Gamma_{I,E}.$
Finally $\tilde\omega_{I,E}=\omega_{a(I,E),b(I,E)}$ (equivalence of the ensembles).
\end{itemize} 

  \end{theo}

\begin{rem}
\label{remap}
$\phantom{a}$

%Before proving the two Theorems above let us remark that,  in this paper,  we 
%dealt  mostly  with the case  where  $E-I$ is fixed.
%The mean field theory in  this case can be handled exactly in the same way.
%In fact, to fix   $ E - I$ for the microcanonical problem,
%is equivalent  to set $a=b$ in the canonical problem, that is to look at the free energy 
%$F(a,b)$ defined above on the line $a=b$.
%The concavity of $F(a, b)$ implies the concavity of $ f(b)=F(b,b)$. 
%This allows us to prove that the solution of the microcanonical problem exists
%and that the entropy
%$S(H)$ ($H=E-I$) is simply given by the Legendre transform of $f(b).$

We finally underline that the function $I \mapsto S(E,I)$ is strictly decreasing
($\partial S/ \partial I = a<0$) so that different radial  solutions of the MFE with the same energy  cannot have
 the same entropy.

\end{rem}

\subsubsection*{Proof of Theorem \ref{theoMFE}}
 
 By the logarithmic
Hardy-Littlewood-Sobolev inequality (see \cite{Beckner}, \cite{CL}), we have
\begin{equation}
\label{LS}
S(\om)-8\pi E(\om)\geq -(1+\log \pi).
\end{equation}
 Note that we also have  the inequality (see \eqref{EI-})
\begin{equation}
\label{EI--}
E(\om)\geq -\frac 1 {8\pi}\log  ( 4 I(\om) ).
\end{equation}
Indeed,  we can write
\begin{equation}
\label{EI-}
E  = -\frac 1{8\pi} \int \log |x-y|^2 \om (x) \om (y)\geq 
-\frac 1{8\pi} \log  \int |x-y|^2 \om (x) \om (y)
\end{equation}
and hence
$$ E \geq  -\frac 1{8\pi} \log \Big( 4 I - 2\Big( \int x \omega \Big) ^2 \Big)
 \geq -\frac 1{8\pi} \log(4I). $$
 Consequently, thanks to \eqref{LS}, \eqref{EI--}, 
we get that
\begin{equation}
\label{F-}
F_{(a,b)} (\omega) \geq   -   {8 \pi - b \over 8 \pi  } \log ( 4 I )  - a I  - ( 1 + \log \pi) 
\end{equation}
and hence, we find that 
$F_{(b,a)}$ is bounded from below. 

Let   $\om_n$  be  a minimizing
sequence in $\Gamma$.  Up to the extraction of a subsequence,
 $\omega_{n}$ converges   in the sense of weak convergence of
measures. Moreover, thanks to \eqref{F-}, we get that $I(\omega_{n})$
 is uniformly  bounded and,   by using again  \eqref{LS}, we also have
 $$( 1 - {b \over 8 \pi }) S(\omega_{n}) \leq  F(\omega_{n}) + { 1 + \log \pi \over 8
 \pi }$$ 
 and hence $S$  is  bounded from above. Therefore the uniform integrability
 given by the bounds on  $S$ and $I$ (which yields
  a bound on $\int \omega | \log \omega | $ thanks to Lemma \ref{lement}) implies  that $\omega_{n}$ converges
  to a nonegative function  $\omega$. Moreover,  the uniform estimate
   on the moment of inertia provides  the tightness of    the sequence $\omega_{n}$  so that we obtain
   $$ \lim_{n} \int \omega_{n} = \int \omega $$
      i.e. $\omega \in \Gamma.$ Next, by lower semi-continuity, we  have
  $$ S(\omega) \leq \lim_{n}S(\omega_{n}), \quad I(\omega) \leq \lim_{n} I(\omega_{n})$$
  and  we claim  that 
  \begin{equation}
  \label{En}
  E(\omega) = \lim_{n} E(\omega_{n}).
  \end{equation}
  This proves that  
  $F_{(a,b)}(\omega) \leq  \lim_{n} F_{(a,b)}(\omega_{n})$  and hence that $\omega $ is a minimizer.
  It remains to prove \eqref{En}.  We write
  $$ E(\omega_{n}) =  -  { 1 \over 4 \pi } \int \log |x- y | \omega_{n}(x) \omega_{n}(y) \, dxdy
   = I(\eps) + J(\eps)$$
   where
   $$ I(\eps) =  -  { 1 \over 4 \pi } \int_{ |x- y |\leq \eps} \log |x- y | \omega_{n}(x) \omega_{n}(y) \, dxdy, \quad  J(\eps)  = -  { 1 \over 4 \pi } \int_{ |x- y |\geq \eps} \log |x- y |
    \omega_{n}(x) \omega_{n}(y) \, dxdy$$
 for every $\eps \in (0, 1 ) $.  By splitting the integration domain in
  $\{ \omega_n(x) \omega_{n}(y) \leq |x- y |^{ -1 } \}$ and its complementary, we easily
   get that
  \begin{eqnarray*}
  I(\eps)  & \leq &   - C \int_{ |x- y | \leq \eps } |x - y |^{-1 } \log |x-y| \, dx dy 
   +  2C\Big( \int  \omega_{n} |\log \omega_{n}| \Big) \sup_{x} \int_{|x- y | \leq \eps}
    \omega_{n}(y)\,  dy \\
    & \leq &  - C \int_{ |x- y | \leq \eps } |x - y |^{-1 } \log |x-y| \, dx dy  + C
    sup_{x} \int_{|x- y | \leq \eps}
    \omega_{n}(y)\,  dy.
    \end{eqnarray*}
     For  the last line, we have used that  the entropy and the moment of inertia are uniformly bounded in $n$ and   
      thanks to  Lemma  \ref{lement}, we also have that 
      $\int \omega_{n} |\log \omega_{n}|$
       is uniformly bounded.  This yields  the uniform integrability :
       $$   \lim_{\eps \rightarrow 0} sup_{x} \int_{|x- y | \leq \eps}
    \omega_{n}(y)\,  dy  = 0.$$
    Consequently, we  get that $\lim_{\eps \rightarrow 0 } I(\eps) = 0$
     uniformly in $n$.
      Finally, by weak convergence, we have that
    $$  \lim_{n}J(\eps) =   - { 1 \over 4 \pi } \int_{ |x- y |\geq \eps} \log |x- y |
    \omega(x) \omega(y) \, dxdy $$
     and since
     $$ \lim_{\eps \rightarrow 0 } - { 1 \over 4 \pi } \int_{ |x- y |\geq \eps} \log |x- y |
    \omega(x) \omega(y) \, dxdy  = E(\omega), $$
    the conclusion follows easily.

\bigskip

By symmetrizing $\om$ (around the origin)  we find that $F_{(a,b)}$ is
decreasing. Indeed $S$ and $I$ are unchanged and $E$ is increasing.
Thus $\om$ must be radially symmetric.
It is not difficult to show that $\om >0$ (otherwise one could find a
better distribution as regards the minimization problem). Hence $\om$
satisfies the MFE. 
\bigskip

Next we show that such a solution is also unique among all the radial
solution to the MFE.
Setting $r=|x|$ and $\psi (r)=\psi (x)$ (by an obvious notational
abuse), we have
\begin{equation}
\frac 1r (r\psi')'=-e^{b\psi +  \frac {a}2 r^2}.
\end{equation}
We are assuming that $Z=2\pi \int_0^\infty dr\, r\, e^{b\psi + \frac {a}2
r^2}=1$, adding, if necessary, a constant to $\psi$.
After the change of variable $t=\log r$, setting $H=b\psi +2t$ we
readily arrive to the following equation:
\begin{equation}
\ddot H=-F(t) e^H
\end{equation}
where
\begin{equation}
F(t)=b e^{a \frac{e^{2t}}2}.
\label{modul}
\end{equation}
We are looking for smooth solutions to eq.n (18) and hence
\begin{equation}
\lim_{t\to-\infty}\dot H=2
\label{asymvel}
\end{equation}
as a consequence of the fact that $\lim _{r\to 0}r\psi'(r)\to 0$.
$H(t)$ behaves as $2t+\chi$ as $t\to-\infty$ and $\chi$ must be
choosen in such a way that $Z=1$. In the new variables:
\begin{equation}
Z=\frac {2\pi}b \int_{-\infty}^\infty dt \, e^H F(t)=\frac {2\pi}b
(\dot H (-\infty)-\dot H (\infty))
\end{equation}
It is convenient to change the time variable by setting $2t \to 2t-\chi$,
so that the problem can be reformulated as
\begin{equation}
\ddot H=-F(t-\frac \chi 2) e^H
\end{equation}
$$
\dot H (-\infty)=2, \, H(t) \approx 2t \quad \hbox {for}Ê
\quad t\to -\infty.
$$
Note that
\begin{equation}
Z(\chi) \to 0 \quad \hbox {for}Ê \quad \chi \to -\infty
\end{equation}
\label{asym1}
and
\begin{equation}
Z(\chi) \to \frac {8\pi}b \quad \hbox {for}Ê \quad \chi \to +\infty.
\label{asym2}
\end{equation}
Eq.n \eqref{asym1} is obvious, while eq.n \eqref{asym2} comes out by 
integrating the Hamiltonian system
\begin{equation}
\ddot H=- e^H
\end{equation}
for which, the energy conservation yields $\dot H (\infty)=-2 $.

Since $\frac {8\pi}b>1$ the value $Z=1$ is certainly taken, at least
once. In order to get uniqueness it remain to show that $\chi \to
Z(\chi)$ is a monotone function, actually it is not decreasing.

Defining
\begin{equation}
G=H +  \frac a2 e^{2(t-\frac \chi 2)},
\end{equation}
we find the following set of non-autonomous equations
\begin{eqnarray}
\ddot H&=&- be^G
\nonumber\\
\ddot G&=&- be^G -4(H-G).
\end{eqnarray}
Note that $H,\dot H$ and $G,\dot G$ satisfy the same condition at
$t=-\infty$.

On the other hand the derivatives $\pa_\chi H=h$ and $\pa_\chi G=g$
satisfy
\begin{eqnarray}
\ddot h&=&- be^G g\nonumber\\
\ddot g&=&(4- be^G)g -4h.
\end{eqnarray}
The conditions at $t\to -\infty$ are vanishing for both $h,\dot h$ and
$g,\dot g$. 

Introducing the energy ${\cal E}=\frac 12 \dot H^2 +b e^G$ we get:
$$
\dot {\cal E}= ba e^G e^{2(t-\frac \chi 2)}\leq 0
$$
and hence
\begin{equation}
be^{G(t)} \leq {\cal E}(t) \leq {\cal E}(-\infty)=2. 
\end{equation}
Therefore $\ddot g \geq 0$ as far as $h\leq 0$ and $\ddot h \leq 0$ as
far as $g\geq 0$. This conditions are indeed verified for $t\approx
-\infty$ so that they are true for all the time. Then $G$ is increasing
as well as
$$
Z=\frac {2\pi}b \int_{-\infty}^\infty dt e^G.
$$
$\Box$
\bigskip

\subsection*{Proof of Theorem \ref{MFEmicro}}

To prove the concavity of $F$ we will prove that,  for any $a_1,b_1,a_2,b_2$ :  
$$F(\frac{a_1+a_2}2,  { b_{1 } + b_{2} \over 2 }) >\frac12 \left(F(a_1, b_{1})+F(a_2, 
 b_{2})\right).$$
Let $a=\frac{a_1+a_2}2,$ and $b=\frac{b_1+b_2}2.$ 
By the linearity of $F_{a, b}(\omega)$ as a function of $a, b$ at fixed $\om$,  we get that
\begin{eqnarray*}F(a, b)&=&F_{a, b}(\omega_{a, b})=
\frac12 F_{a_1, b_{1}}(\omega_{a, b})+\frac12 F_{a_2, b_{2}}(\omega_{a, b})\\
&>&\frac12 F_{a_1, b_{1}}(\omega_{a_1, b_{1}})+\frac12 F_{a_2, b_{2}}
(\omega_{a_2, b_{2}})
=\frac12 (F(a_1, b_{1})+F(a_2, b_{2}) ) ,\end{eqnarray*}
where we used the fact that $\omega_{a_1, b_{1}}$ and $\omega_{a_2, b_{2}}$ 
are the minimizer for 
$F_{a_1, b_{1}}$ and $F_{a_2, b_{2}}$ respectively.

The smoothness of $F$ comes from the fact that the solution of the 
canonical variational principle  depends smoothly upon $a,b.$ 
By taking the derivative of $F$ with respect to $a$ we get 
$$\frac{\partial F}{\partial a}=\frac{\partial F_{a,b}(\omega_{a,b})}
{\partial a}=- I(\omega_{a, b}).$$
Here we have used the fact that the derivative of $F$ with respect to 
$\omega$ evaluated 
in $\omega_{a, b}$ vanishes, and the fact that the derivative of $F_{a,b}$ 
with respect to the parameter $a$
is given by $I.$
In the same way we get $\partial F(a, b)/\partial b=-E(\omega_{a,b}).$

Finally, the concavity of  $F_{a, b}$ implies  that
$\partial I/\partial a>0,$ and that $\partial E/\partial b>0$
 again  with the notation  $I(a,b)= I(\om_{a,b}), \, E(a,b) = E(\om_{a,b})$.

Now it remains to prove iii). Again,  the concavity of $F(a, b)$ 
implies the existence
 of  the  convex function  $S^*(I,E)$ defined in \eqref{Lgdr}.  
 Now we want to prove that $S(I,E)=S^*(I,E).$
First of all let us notice that
$$S^*(I,E)= \sup_{a, b}\Big( F(a, b) + b E +a I\Big)=F(\bar a,\bar b) +
\bar b E+\bar a I,$$
where $\bar a,\bar{b}$ is the unique maximum point for $S^*(I,E).$
Therefore, for any $a,b$ 
$$S^*(I,E)\geq F(a, b)+b E +a I=S(\omega_{a, b})-b(E(\omega_{a, b})-E)
-a(I(\omega_{a,b})-I).$$
Now, since $F$ is concave and smooth, we know that for any $I,E,$ there exists unique 
$a,b$ such that $I(\omega_{a, b})=I,$ and $E(\omega_{a,b})=E.$
By choosing $a, b$ in this way in the previous equation we get
\begin{equation}S^*(I,E)\geq S(\omega_{a, b})\geq S(I,E).
\label{Sabove}\end{equation}

On the other hand, let  $\omega_k:k=1,2,...$ be a minimizing sequence for $S(I,E),$ and
$\omega$ a limit point for it. By lower semicontinuity of $S$ we know that
$S(\omega)\leq S(I,E).$ 
Therefore, for any $a,b,$
\begin{eqnarray}
S(I,E)\geq S(\omega)&=&S(\omega)-b E(\omega)-a I(\omega)+b E(\omega)+a I(\omega)\nonumber\\
&\geq& F(a, b)+b E(\omega)+a I(\omega)\geq F(a, b)+b E +a I,\label{vicemic}
\end{eqnarray}
where we have used the continuity of $E$ from which $E(\omega)=E,$ 
and the lower semicontinuity
of $I,$ from which $I(\omega)\leq I.$

Since  $a,b$ are  arbitrary in \eqref{vicemic},  we get 
$S(I,E)\geq S^*(I,E),$.  Since we have already proven   \eqref{Sabove} 
this yields  $S(I,E)=S^*(I,E).$

Finally let us notice that 
$$S(I,E)=S^*(I,E)=S(\omega_{\bar a,\bar{b}}),$$
where $\bar{a},\bar{b}$ is the unique minimun point for $S^*(I,E),$ 
where the relation between $a,b$ and $I,E$ is smooth and bijective 
(equivalence of the ensembles).

\bigskip

\bibliographystyle{acm}
\bibliography{mcan}

\end{document}